\documentclass{amsart}

\author{Samuel Boissi{\`e}re \and Alessandra Sarti}

\title{Counting lines on surfaces}

\date{5th June 2006}

\address{Samuel Boissi{\`e}re, Laboratoire J.A.Dieudonn\'e UMR CNRS 6621,
         Universit\'e de Nice Sophia-Antipolis, Parc Valrose, 06108 Nice}

\email{sb@math.unice.fr}

\urladdr{http://math.unice.fr/$\sim$sb/}

\address{Alessandra Sarti, Fachbereich f{\"u}r Mathematik, Johannes Gutenberg-Universit{\"a}t,
55099 Mainz, Germany}
\curraddr{Dipartimento di Matematica, Universit\`a di Milano, via Saldini 50, 20133 Milano, Italy}
\thanks{The second author is supported by DFG Research Grant SA 1380/1-2.}

\email{sarti@mathematik.uni-mainz.de, sarti@mat.unimi.it}

\urladdr{http://www.mathematik.uni-mainz.de/$\sim$sarti}

\usepackage[english]{babel}
\usepackage{amsmath,amsfonts,amssymb}
\usepackage{bbm,dsfont}
\usepackage{graphicx}

\DeclareMathOperator{\Span}{Span} \DeclareMathOperator{\Pic}{Pic}
\DeclareMathOperator{\PGL}{PGL} \DeclareMathOperator{\diag}{diag}
\DeclareMathOperator{\id}{id} \DeclareMathOperator{\SO}{SO}
\DeclareMathOperator{\SU}{SU}

\newcommand{\ie}{\textit{i.e. }}
\newcommand{\eg}{\textit{e.g. }}
\newcommand{\apriori}{\textit{a priori }}
\newcommand{\loccit}{\textit{loc.cit. }}

\newcommand{\ii}{\mathrm i}
\newcommand{\jj}{\omega}

\newcommand{\cC}{\mathcal{C}}

\newcommand{\cH}{\mathcal{H}}
\newcommand{\cI}{\mathcal{I}}
\newcommand{\cO}{\mathcal{O}}
\newcommand{\cR}{\mathcal{R}}
\newcommand{\cS}{\mathcal{S}}
\newcommand{\cT}{\mathcal{T}}
\newcommand{\cU}{\mathcal{U}}

\newcommand{\fA}{\mathfrak{A}}

\newcommand{\fS}{\mathfrak{S}}

\newcommand{\IC}{\mathds{C}}
\newcommand{\IQ}{\mathds{Q}}
\newcommand{\IR}{\mathds{R}}
\newcommand{\IZ}{\mathds{Z}}

\newcommand{\IG}{\mathbb{G}}
\newcommand{\II}{\mathbb{I}}
\newcommand{\IP}{\mathbb{P}}

\newtheorem{theorem}{Theorem}[section]

\newtheorem{proposition}[theorem]{Proposition}
\newtheorem{corollary}[theorem]{Corollary}
\newtheorem{remark}[theorem]{Remark}

\begin{document}

\begin{abstract}
This paper deals with surfaces with many lines. It is well-known that a cubic
contains $27$ of them and that the maximal number for a quartic is $64$. In
higher degree the question remains open. Here we study classical and new
constructions of surfaces with high number of lines. We obtain in particular a
symmetric octic with $352$ lines.
\end{abstract}



\maketitle

\vspace{-1cm}

\begin{center}
\begin{figure}[h!]
\includegraphics[width=4cm, height=4cm]{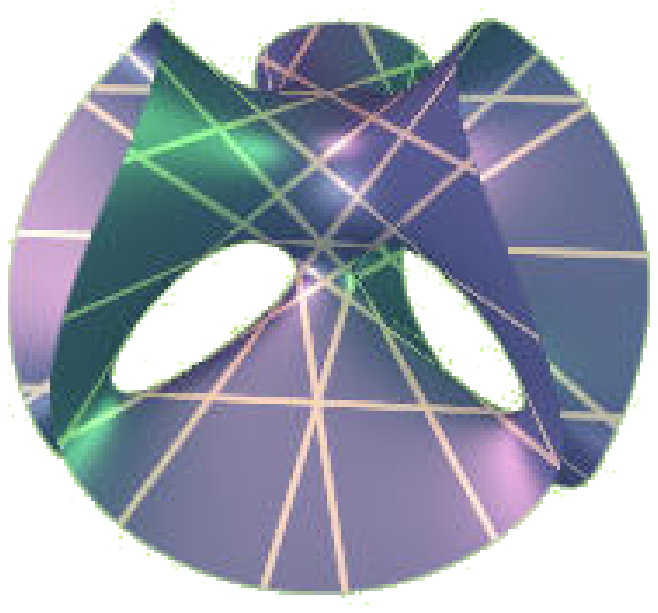}
\end{figure}
\vspace{-0.4cm} {Cubic surface with $27$ lines\footnote{{\tt
http://enriques.mathematik.uni-mainz.de/surf/logo.jpg}}}
\end{center}

\section{Introduction}

Motivation for this paper is the article  in 1943 of Segre \cite{Segre43} which
studies the following classical problem: What is the maximum number of lines a
surface of degree $d$ in $\IP_3$ can have? Segre answers this question for
$d=4$ by using some nice geometry, showing that it is exactly $64$. For the
degree three it is a classical result that each smooth cubic in $\IP_3$
contains $27$ lines, but for $d\geq 5$ this number is still not known. In this
case, Segre shows in \loccit that the maximal number is less or equal to
$(d-2)(11d-6)$ but this bound is far from beeing sharp. Indeed, already in
degree four it gives $76$ lines which is not optimal. So from one hand one can
try to improve the upper bound for the number of lines $\ell(d)$ a surface of
degree $d$ in $\IP_3$ can have, on the other hand it is interesting to
construct surfaces with as many lines as possible to give a lower bound for
$\ell(d)$.

It is notoriously difficult to construct examples of surfaces with many lines.
Good examples so far are the surfaces of the kind
$F(x,y,z,t)=\phi(x,y)-\psi(z,t)=0$ where $\phi$ and $\psi$ are homogeneous
polynomials of degree $d$. Segre in \cite{Segre47} studies the case of $\deg
F=4$ showing that in this case the possible numbers of lines are $16,32,48,64$.
He finds these numbers by studying the automorphisms of $\IP_1$ between the two
sets of four points $\phi=0$ and $\psi=0$. Caporaso-Harris-Mazur in \cite{CHM},
 by using similar methods as Segre, then study the maximal number of lines $N_d$ 
on such surfaces in any degree $d$ showing that $N_d\geq 3d^2$ for each $d$ and $N_{4}\geq 64$, $N_6\geq 180$,
$N_8\geq 256$, $N_{12}\geq 864$, $N_{20}\geq 1600$. In this paper we show the exactness of these results.  First we note
that it is enough to consider surfaces of the kind $\phi(x,y)-\phi(z,t)=0$ and
by a careful analysis of the automorphisms of the set of points $\phi=0$ on
$\IP_1$ we can list all the possible numbers of lines on surfaces of this kind
for all $d$ and then we prove:

\noindent{\bf Proposition \ref{prop:aziz} }{\it The maximal numbers of lines on $F=0$ are:
\begin{itemize}
\item $N_d=3d^2$ for $d\geq 3$, $d\not=4,6,8,12,20$;

\item $N_4=64$, $N_6=180$, $N_8=256$, $N_{12}=864$, $N_{20}=1600$.
\end{itemize}}

It is well-known that the Fermat surfaces $(x^d-y^d)-(z^d-t^d)=0$ have $3d^2$
lines. Our proof provides a method to write equations of surfaces
$\phi(x,y)-\psi(z,t)=0$ with each possible number of lines. In particular, our
proposition shows that it is not possible, with these surfaces, to obtain
better examples and a better lower bound for $\ell(d)$. So, in order to find
better examples, one has to use new methods. In this paper we explore the
following kinds of surfaces:
\begin{itemize}
\item $d$-covering of the plane $\IP_2$ branched over a curve of degree $d$.

\item Symmetric surfaces in $\IP_3$.
\end{itemize}
We show that the first method cannot give more than $3d^2$ lines (Proposition
\ref{prop:covering}).

The second method is based on the following idea: if a surface has many
automorphisms (many symmetries) then possibly it contains many orbits of lines.
This idea was used successfully in the study of surfaces with many nodes. In
this paper we find a $G_8$-invariant octic with $352$ lines, where $G_8\subset
\PGL(3,\IC)$ has order $576$ (Proposition \ref{prop:record}). This shows
$\ell(8)\geq 352$, improving the previous bound of $256$.

As stated before, one can also try to improve the upper bound for $\ell(d)$.
Following the idea of Segre \cite{Segre43} and imposing some extra conditions
on the lines on a surface, we can find the bound $d(7d-12)$ which surprisingly
agrees with the maximal examples in degrees $4,6,8,12$ (Section
\ref{s:record}).

Finally a related problem to this is to determine the maximal number $m(d)$ of
skew-lines a surface of degree $d$ in $\IP_3$ can have. It is well-known that
$m(3)=6$ and $m(4)=16$. For $d\geq 5$, this value is not known. An upper bound
$m(d)\leq 2d(d-2)$ is given by Miyaoka in \cite{Mi}, which is sharp for
$d=3,4$. There are results of Rams \cite{Rams02, Rams05} giving examples of
surfaces with $d(d-2)+2$ skew-lines ($d\geq 5$) and with $19$ skew-lines for
$d=5$. In Proposition \ref{essla1} we improve his examples for $d\geq 7$ and
$\gcd(d,d-2)=1$ to $d(d-2)+4$.

The paper is organized as follows. In Section \ref{ch:general} we give an
overview of known results. In Sections \ref{ch:segre} and \ref{ch:geemen} we
describe completely the surfaces of the kind $\phi(x,y)-\psi(z,t)=0$ and the
$d$-coverings of the plane $t^d=f(x,y,z)$. Section \ref{ch:symmetric} is
devoted to the investigation of symmetric surfaces, and in particular of an
octic with $352$ lines. In Section \ref{s:record} we present the uniform bound
$d(7d-12)$ and Section \ref{s:aziz} is an application to the problem of the
number of rational points on curves. Finally, Section \ref{s:skew} deals with
the skew-lines: we give an overview of known results and some new examples.

\bigskip

{\it Acknowledgements.} We thank Duco van Straten for suggesting us this nice problem and for interesting discussions.

\section{General results}
\label{ch:general}

Our objective is to investigate the number of lines
contained in a smooth surface in $\IP_3$. We first recall classical
results: the generic situation and the bound of Segre.

\subsection{Generic situation}\text{}
\label{s:generic}

It is a well-known fact that each smooth quadric surface in $\IP_3$ contains an
infinite number of lines and each smooth cubic surface in $\IP_3$ contains
exactly $27$ lines. What happens for surfaces of higher degree? Generically:

\begin{proposition}
A generic smooth surface of degree $d\geq 4$ in $\IP_3$ contains no line.
\end{proposition}

We briefly recall the proof, following \cite{AK,BV}.

\begin{proof}\text{}
Let $V$ be the vector space of degree $d$ homogeneous polynomials in the coordinates $x,y,z,t$ and $G$ be the Grassmannian of $2$-planes in $\IC^4$. Consider the incidence variety $F:=\{(L,f)\subset G\times V\,|\,f_{|L}\equiv 0\}$ with its projections $p:F\rightarrow G$ and $q:F\rightarrow V$.

$\bullet$ Let $L\in G$ and assume that $L$ is generated by the vectors $(1,0,0,0)$ and $(0,1,0,0)$ in $\IC^4$. Consider the affine neighbourhood of $L$ in $G$:
$$
\cU:= \Span\{(1,0,a,b),(0,1,c,d)\}
$$
where $a,b,c,d$ are local coordinates. If $f\in p^{-1}(\cU)$, then
$$
f\left(\lambda(1,0,a,b)+\mu(0,1,c,d)\right)=0\qquad \forall \lambda,\mu\in \IC,
$$
and denoting $f=\sum\limits_{i+j+k+l=d}a_{i,j,k,l}x^iy^jz^kt^l$, the equation
$$
\sum\limits_{i+j+k+l=d}a_{i,j,k,l}\lambda^i\mu^j(\lambda a+\mu c)^k(\lambda b+\mu d)^l=0 \qquad\forall \lambda,\mu\in \IC
$$
gives $d+1$ linear equations in the coordinates $(a_{i,j,k,l})$ of $f\in V$ whose rank at $a=b=c=d=0$ is $d+1$: hence locally in a neighbouhood of $L$, the system has rank $d+1$ so $p$ is a locally trivial bundle of rank: $\dim V -(d+1)$.

$\bullet$ Let $X$ be a surface of degree $d$ in $\IP_3$, given by a polynomial $f\in V$. Then the \emph{Fano scheme} parametrizing the lines contained in $X$ is $F(X):=~p\left(q^{-1}(f)\right)$.

$\bullet$ Consider the map $q:F\rightarrow V$. Since:
$$\dim F=\dim V-(d+1)+\dim G=\dim V -(d-3),
$$
for $d\geq 4$ one has $\dim F<\dim V$ hence the map $q$ is not dominant. This means that the generic fibre of $q$ is empty. Otherwise stated, $F(X)$ is empty for $X$ generic.
\end{proof}

We shall see in the next section that the number of lines a smooth surface of
degree $d\geq 4$ can have is always finite, and bounded. This leads to the
problem of finding surfaces with an optimal number of lines.

\subsection{Upper bound for lines}\text{}
\label{s:upperline}

The best upper bound known so far for the number of lines on a smooth
surface of degree $d\geq 4$ in $\IP_3$ is given by Segre:

\begin{theorem}[Segre \cite{Segre43}]\text{}
\begin{itemize}
\item The number of lines lying on a smooth surface of degree $d\geq 4$
does not exceed $(d-2)(11d-6)$.

\item The maximum number of lines lying on a quartic
surface is exactly $64$.
\end{itemize}
\end{theorem}

This bound is effective for $d=4$ (see for instance maximal examples in Section
\ref{s:segregeneral}) but for $d\geq 5$ it is believed that it could be
improved. For instance, already for $d=4$ the uniform bound $(d-2)(11d-6)$ is
too big. The next sections are devoted to the study of some families of
surfaces with particular properties, containing many lines.

\section{Surfaces of the kind $\phi(x,y)=\psi(z,t)$}
\label{ch:segre}

We consider a surface $\cS$ given by an equation of the kind:
$$
F(x,y,z,t):=\phi(x,y)-\psi(z,t)
$$
for two homogeneous polynomials $\phi,\psi$ of degree $d$. Segre gave a
complete description of the possible and maximal numbers of lines in the case
$d=4$ (\cite[\S VIII]{Segre47}). We generalize the method to all degrees: we
treat in details the configuration of lines, give a description of all
possible numbers, and conclude with the maximal numbers of lines for
such surfaces.

\subsection{Configuration of the lines}\text{}
\label{s:segregeneral}

Let $Z(\phi)$, resp. $Z(\psi)$ denote the set of zeros of $\phi(x,y)$, resp.
$\psi(z,t)$ in $\IP_1$.

\begin{theorem}\label{teoremone}
Let $F(x,y,z,t)=\phi(x,y)-\psi(z,t)$ be the equation of a smooth surface $\cS$
of degree $d$ in $\IP_3$. The number $N_d$ of lines on $\cS$ is exactly:
\begin{eqnarray*}
N_d=d(d+\alpha_d)
\end{eqnarray*}
where $\alpha_d$ is the order of the group of isomorphisms of $\IP_1$ mapping
$Z(\phi)$ to $Z(\psi)$.
\end{theorem}

\begin{proof}\text{}

$\bullet$ Let $L$ be the line $z=t=0$ and $L'$ be the line $x=y=0$. Then
$\cS\cap L=Z(\phi)$ and $\cS\cap L'=Z(\psi)$. Since the surface $\cS$ is
smooth, the homogeneous polynomials $\phi$ and $\psi$ have simple zeros.
Indeed, for example in the case of the polynomial $\phi$, if $[a:b]\in \IP_1$
is such that $\phi$ can be factorized by $(bx-ay)^2$, then
$\partial_x\phi(a,b)=\partial_y\phi(a,b)=0$ and the point $[a:b:0:0]$ is a
singular point of $\cS$ (the inverse also holds: if both $\phi$ and $\psi$ have
only simple zeros, then $\cS$ is smooth). Set $Z(\phi):=\{P_1,\ldots,P_d\}$ and
$Z(\psi):=\{P'_1,\ldots,P'_d\}$.

$\bullet$ Each line $L_{i,j}$ joining a $P_i$ to a $P'_j$ is contained in
$\cS$: if $P_i=[x_i:y_i:0:0]$ and $P'_j=[0:0:z'_j:t'_j]$ the line joining
them consists in points $[\lambda x_i:~\lambda y_i:~\mu z'_j:~\mu
t'_j]$, $\lambda,\mu\in \IC$, which are all contained in the surface, by
homogeneity of the polynomials $\phi$ and $\psi$. This gives $d^2$ lines.

$\bullet$ Each line contained in $\cS$ and intersecting $L$ and $L'$ is one of
the previous lines. Indeed, if $D$ is such a line, set $D\cap L=\{[a:b:0:0]\}$
and $D\cap L'=\{[0:0:c:d]\}$. Then $F(a,b,0,0)=\phi(a,b)=0$ so $[a:b:0:0]$ is
one of the points $P_i$ and similarly $[0:0:c:d]$ is one $P'_j$.

$\bullet$ Let $D$ be a line contained in $\cS$ and not intersecting $L$. Then
$D$ does not intersect $L'$ (and {\it vice-versa}). Indeed, an equation of such
a line $D$ is given by two independent equations:
\begin{equation*}
\left\{
\begin{aligned}
ax+by+cz+dt &=0\\
a'x+b'y+c'z+d't &=0
\end{aligned}
\right.
\end{equation*}
Since $D$ does not intersect $L$, the system
\begin{equation*}
\left\{
\begin{aligned}
ax+by &=0\\
a'x+b'y &=0
\end{aligned}
\right.
\end{equation*}
has rank two, so we can rewrite the equations of $D$ as the following
independent equations:
\begin{equation*}
\left\{
\begin{aligned}
x&=\alpha z+\beta t\\
y&=\gamma z+\delta t
\end{aligned}
\right.
\end{equation*}
Then $D$ does not intersect $L'$ otherwise the matrix $\left(
\begin{array}{cc}
  \alpha & \beta \\
  \gamma & \delta
\end{array}
\right)$ would have rank one.

$\bullet$ Therefore, the equations of the line $D$ define a linear isomorphism
between the lines $L'$ and $L$ inducing a bijection between
 $Z(\psi)$ and $Z(\phi)$. Indeed, seting $P'_j=[0:0:c:d]$, then
$a:=\alpha c+\beta d$ and $b:=\gamma c+\delta d$ have the property that
$[a:b:c:d]\in D\subset \cS$ so $\phi(a,b)=F(a,b,c,d)+\psi(c,d)=0$
hence $[a:b:0:0]$ is a zero of $\phi$.

$\bullet$ Conversely, let $\sigma:L'\rightarrow L$ be an isomorphism mapping the points $P'_j$ to the points $P_i$,
and $\left(\begin{array}{cc} \alpha& \beta\\ \gamma & \delta \end{array}\right)$ a matrix defining $\sigma$. Consider the smooth quadric $Q_\sigma : x(\gamma z+\delta t)-y(\alpha z+\beta t)=0$.
Its first ruling is the family of lines $(p,\sigma(p))$ for $p\in L'$. For $p=[c:d]$, these lines are given by the equations
$$
\mathrm{I}_{[c:d]}:\left\{
\begin{aligned}
(\gamma c +\delta d)x-(\alpha c+\beta d)y&=0\\
dz-ct&=0
\end{aligned}
\right.
$$
Its second ruling consists in the family of lines of equations
$$
\II_{[a:b]}:\left\{
\begin{aligned}
ax-b(\alpha z+\beta t)&=0\\
ay-b(\gamma z+\delta t)&=0
\end{aligned}
\right.
$$
for $[a:b]\in \IP_1$. To this ruling belong the lines $L$ ($[a:b]=[0:1]$), $L'$
($[a:b]=[1:0]$) and $D$ ($[a:b]=[1:1]$). It is not true \apriori that
\emph{this} $D$ is contained in $\cS$, since the matrix $\sigma$ is defined up
to a scalar factor.

In each ruling, the lines are disjoint to each other, and each line of one
ruling intersects each line of the other ruling. Since the intersection
$\cS\cap Q$ contains exactly the $d$ different lines $(P'_j,\sigma(P'_j))$ of
the first ruling, it contains also $d$ lines of the second ruling: Consider a
line in the first ruling not contained in $\cS$, then it intersects $\cS$ in
$d$ points, and through each of this points is attached a line of the second
ruling, which also intersects the $d$ lines of the first ruling contained in
$\cS$, so these lines of the second ruling intersect $\cS$ at $d+1$ points, so
are contained in $\cS$. But it is not clear \apriori with our argument that
these lines in the second ruling are different. Denote by $\cU_d$ the group of
$d$-th roots of the unit. The group $\cU_d\times \cU_d$ acts on $\IP_3$ by
$(\xi,\eta)\cdot [x:y:z:t]=[\xi x:\xi y:\eta z:\eta t]$, leaving the surface
$\cS$ globally invariant since the polynomials $\phi$ and $\psi$ are
homogeneous of degree $d$. Observe that the lines of the first ruling are
invariant for the action, but for the second ruling, $(\xi,\eta)\cdot
\II_{[a:b]}=\II_{[\xi^{-1}a:\eta^{-1}b]}$ so each line of the second ruling
produces a length $d$ orbit through the action. Since the surface $\cS$
contains at least one line of the second ruling, it contains the whole orbit,
this gives us $d$ different lines.

Therefore, each isomorphism $\sigma:L'\rightarrow L$ mapping $Z(\psi)$ to
$Z(\phi)$ gives $d$ lines, and there are no other lines. Furthermore, for two
different isomorphisms, the corresponding lines are different since the matrix
defining the isomorphims are not proportional.

$\bullet$ Denote by $\alpha_d$ the number of isomorphims $\sigma:L'\rightarrow
L$ mapping $Z(\psi)$ to $Z(\phi)$. The preceding discussion shows that the
exact number of lines contained in the surface $\cS$ is:
$$
N_d=d^2+\alpha_d d.
$$
\end{proof}

\begin{remark}
In the proof of \cite[Lemma 5.1]{CHM}, Caporaso-Harris-Mazur proved with a
similar argument that the number of lines is at least $d(d+\alpha_d)$ and
described some special values. Our argument includes the exactness. In the next
subsections we give a full description of the possible values of $\alpha_d$, in
particular its maximal values for each $d$.
\end{remark}

\subsection{The possible numbers of lines}\text{}
\label{ss:possiblesegre}

Now we want to find the possible and maximal values of $N_d$, or equivalently
$\alpha_d$. If there is at least one isomorphism $\sigma$ (see the proof
above), then by composing by $\sigma^{-1}$ we are lead to the problem of
determining the possible numbers of automorphisms of $\IP_1$ (or
\emph{projectivities}) acting on a given set of $d$ points on $\IP_1$. Since a
projectivity is defined by its value on three points, we have always
$\alpha_3=6$, and for $d\geq 4$ there is only a finite number of such
isomorphisms, depending on the relative position of the points, encoded in
their cross-ratios. The case $d=4$ was studied by Segre \cite{Segre47} with
this point of view. We give a different argument for the general case. The set
$\Gamma_d$ of isomorphims of $\IP_1$ acting on  $d$ points defines a finite
group of automorphisms of $\IP_1$. First recall the classical classification:

\bigskip

{\bf Polyhedral groups.} There are five types of finite subgroups of
$\SO(3,\IR)$, or equivalently of $\PGL(2,\IC)$, called \emph{polyhedral
groups}:
\begin{itemize}
\item the cyclic groups $C_k\cong \IZ/k\IZ$ of order $k\geq 2$, isomorphic to
the group of isometries of a regular polygon with $k$ vertices in the plane;

\item the dihedral groups $D_k\cong \IZ/k\IZ \rtimes \IZ/2\IZ$ of order $2k$,
$k\geq 2$, isomorphic to the group of isometries of regular polygon with $k$
vertices in the space;

\item the group $\cT$ of positive isometries of a regular tetrahedra,
isomorphic to the alternate group $\fA_4$ of order twelve;

\item the group $\cO$ of positive isometries of a regular octahedra or a cube,
isomorphic to the symmetric group $\fS_4$ of order $24$;

\item the group $\cI$ of positive isometries of a regular icosahedra or a
regular dodecahedra, isomorphic to the alternate group $\fA_5$ of order $60$.
\end{itemize}
In the sequel, we shall describe generators of these groups and
their orbits on $\IP_1$, in order to get explicit constructions of surfaces.

\bigskip

We now proceed to the description of all possible groups of isomorphisms ($d\geq 4$):
\begin{enumerate}
\item {\it $\Gamma_d=\{\id\}$.} This is not possible for $d=4$ since there are
always at least four automorphisms of a set of four points in $\IP_1$ (their cross-ratio
takes generically six different values under permutation).

\item {\it $\Gamma_d$ is a cyclic group:} $\Gamma_d\cong \IZ/k\IZ$ ($k\geq
2$) with generator $\sigma(t)=\xi t$ where $\xi$ is a primitive $k$-th root of
the unit. The action of $\sigma$ on $\IP_1$ has two fix points $\{0,\infty\}$ and all other points generate a length $k$ orbit. So, depending
whether the fix points are in the given set of $d$ points or not we have the decomposition:
$$
d=\alpha+\beta k
$$
with $\alpha\in\{0,1,2\}$ and $\beta\geq 1$:
\begin{itemize}
\item $\alpha=0$. The points are\footnote{Here and in the sequel, the $\mu_i$'s
are assumed to be generic: they are distinct and in particular they are not the
$\beta$-th roots of the unit and their $k$-powers $\lambda_i:=\mu_i^k$ are
distincts.}:
$$
\{\mu_1,\mu_1\xi,\ldots,\mu_1\xi^{k-1}\},\ldots,\{\mu_\beta,\mu_\beta\xi,\ldots,\mu_\beta\xi^{k-1}\}.
$$
This forces $\beta\geq 3$ since: if $\beta=1$ or $\beta=2$ then
$t\mapsto 1/t$ or $t\mapsto \mu_2/(\mu_1t)$ generate a dihedral group. For $\beta\geq 3$ there are no
other isomorphisms.

\item $\alpha=1$. The points are:
$$
\{0\},\{\mu_1,\mu_1\xi,\ldots,\mu_1\xi^{k-1}\},\ldots,\{\mu_\beta,\mu_\beta\xi,\ldots,\mu_\beta\xi^{k-1}\}.
$$
There is no other isomorphism whenever $d=1+\beta k\geq 5$. For $k=3$ and
$\beta=1$ there are other isomorphism (a tetrahedral group).

\item $\alpha=2$. The points are:
$$
\{0,\infty\},\{\mu_1,\mu_1\xi,\ldots,\mu_1\xi^{k-1}\},\ldots,\{\mu_\beta,\mu_\beta\xi,\ldots,\mu_\beta\xi^{k-1}\}.
$$
As before, this forces $\beta\geq 3$.
\end{itemize}

To summarize, for the group $\Gamma_d$ be a cyclic group $\IZ/k\IZ$ ($d\geq
4$, $k\geq 2$):
\begin{itemize}
\item $d=\beta k$, $\beta\geq 3$, \eg
$\phi(x,y)=\prod\limits_{i=1}^\beta(x^k-\lambda_iy^k)$;

\item $d=1+\beta k\geq 5$, $\beta\geq 1$ if $k=3$, \eg
$\phi(x,y)=x\prod\limits_{i=1}^\beta(x^k-\lambda_iy^k)$;

\item $d=2+\beta k$, $\beta\geq 3$ \eg
$\phi(x,y)=xy\prod\limits_{i=1}^\beta(x^k-\lambda_iy^k)$.
\end{itemize}

\item {\it $\Gamma_d$ is a dihedral group:} $\Gamma_d\cong \IZ/k\IZ\rtimes
\IZ/2\IZ$ ($k\geq 2$) with generators $\sigma(t)=\xi t$ and $s(t)=1/t$  where
$\xi$ is a primitive $k$-th root of the unit. The action of the
dihedral group on $\IP_1$ has one length 2 orbit $\{0,\infty\}$ and one length $k$ orbit generated by $1$. So we have the decomposition:
$$
d=2\alpha+\beta k+\gamma 2k
$$
with $\alpha,\beta\in\{0,1\}$, $\gamma\geq 0$:
\begin{itemize}
\item $\gamma=0$, $\alpha=0$ and $\beta=1$. The points are:
$$
\{1,\xi,\ldots,\xi^{k-1}\}
$$
Then $d=k$ and $\phi(x,y)=x^k-y^k$. This gives the Fermat surface.

\item $\gamma=0$, $\alpha=1$ and $\beta=1$. The points are:
$$
\{0,\infty\},\{1,\xi,\ldots,\xi^{k-1}\}
$$
This forces $k\neq 2,4$: if $k=2$, the configuration is isomorphic to the preceding case (with $2k$) and contains more isomorphims, and if $k=4$ there are other isomorphisms generating an octahedral group. Then $d=2+k$ and $\phi(x,y)=xy(x^k-y^k)$.

\item $\gamma\neq 0$. Then $d\in\{2k\gamma,2+2k\gamma, k+2k\gamma,2+k+2k\gamma
\}$ and $\phi$ contains, besides the factors given in the preceding cases,
$\gamma$ factors of the kind $(x^k-\lambda y^k)(x^k-\frac{1}{\lambda}y^k)$.
\end{itemize}

\item {\it $\Gamma_d$ is a tetrahedral group $\cT$.} The group $\cT$ is
generated by:
$$
\sigma(t)=\jj t,\quad s(t)=\frac{1-t}{1+2t}
$$
acting on the set $\{0,1,\jj,\jj^2\}$ where $\jj$ is a primitive third root of
the unit. The action of $\cT$ on $\IP_1$ has two length four orbits:
$$
\{0,1,\jj,\jj^2\},\{\infty,\frac{-1}{2},\frac{-1}{2}\jj,\frac{-1}{2}\jj^2\}
$$
and one length six orbit generated by the fix point $w=\frac{-1-\sqrt{3}}{2}$
of\footnote{The second fix point $w'=\frac{-1+\sqrt{3}}{2}$ belongs to the same
orbit since $w=\sigma^2s\sigma(w')$.} $s$. These are all the orbits of lengths
four or six since the conjugacy classes in $\cT$ are generated by
$\id,s,\sigma,\sigma^2$. So we have the decomposition:
$$
d=4\alpha+6\beta+12\gamma
$$
with $\alpha\in\{0,1,2\}, \beta\in\{0,1\}, \gamma\geq 0$:
\begin{itemize}
\item $\gamma=0$, $\beta=0$ and $\alpha=1$: the group of isomorphisms is $\cT$.

\item $\gamma=0$, $\beta=0$ and $\alpha=2$: the group of isomorphisms would be
$\cO$ since $t\mapsto -1/(2t)$ interchanges the two length four orbits.

\item $\gamma=0$, $\beta=1$ and $\alpha=0$: the group of isomorphisms would be
$\cO$ since the length six orbit is stabilized by $t\mapsto  -1/(2t)$.

\item $\gamma=0$, $\beta=1$ and $\alpha=1$: the group of isomorphisms is $\cT$,
because it is not contained in any dihedral group and the groups $\cO$ or $\cI$
have no length four or ten orbit.

\item $\gamma=0$, $\beta=1$ and $\alpha=2$: as before the group of isomorphisms is
$\cO$.

\item For $\gamma\neq 0$, in general the group of isomorphisms is $\cT$ but for
special points this could be $\cO$ or $\cI$.
\end{itemize}

For example, for the tetrahedral group consider $\phi(x,y)=x(x^3-y^3)$.

\item {\it $\Gamma_d$ is an octahedral group $\cO$.} The group $\cO$ is
generated by:
$$
\sigma(t)=\ii t,\quad s(t)=\frac{1}{t},\quad a(t)=\frac{t+\ii}{t-\ii}
$$
acting on the set $\{0,\infty,1,\ii,-1,-\ii\}$. The action of $\cO$ on $\IP_1$
has one length six orbit, one length eight orbit generated by the fix point
$w=\frac{1+\ii-\sqrt{3}-\ii\sqrt{3}}{2}$ of\footnote{The second fix point
$w'=\frac{1+\ii+\sqrt{3}+\ii\sqrt{3}}{2}$ belongs to the same orbit since $w'=s
a \sigma s(w)$.} $a$, and one length twelve orbit generated by the fix point
$z=-1+\sqrt{2}$ of the isomorphism\footnote{The second fix point
$z'=-1-\sqrt{2}$ belongs to the same orbit since $z'=\sigma r \sigma a(z)$.}
$r(t)=\frac{1-t}{1+t}$. These are all orbits of lengths six, eight or twelve
since the conjugacy classes in $\cO$ are generated by $\id,s,\sigma,a,r$. So we
have the decomposition:
$$
d=6\alpha+8\beta+12\gamma+24\delta
$$
with $\alpha,\beta,\gamma\in\{0,1\}$, $\delta\geq 0$. Since the group $\cO$ is
not contained in $\cI$ nor in any dihedral group, all choices of
$\alpha,\beta,\gamma,\delta$ are possible to get $\Gamma_d\cong \cO$.

\item {\it $\Gamma_d$ is a icosahedral group $\cI$.} The group $\cI$ is
generated by:
$$
p_5(t):=\frac{\tau
t+\tau-1+\ii}{(-\tau+1+\ii)t+\tau},\quad q_1(t):=-t,\quad q_2(t):=-\frac{1}{t}
$$
where $\tau:=\frac{1+\sqrt{5}}{2}$. The only length twelve orbit is generated
by a fix point of $p_5$, the length $20$ orbit is generated by a fix point of
$p_5^2q_2$ (which has order three) and the length $30$ orbit is generated by a
fix point of $q_1$. Since the conjugacy classes in $\cI$ are generated by
$\id,p_5,p_5^2,p_5^2q_2,q_1$ there are no other orbits. So we have the
decomposition:
$$
d=12\alpha+20\beta+30\gamma+60\delta
$$
with $\alpha,\beta,\gamma\in\{0,1\}$, $\delta\geq 0$. All choices give $\Gamma_d\cong \cI$.
\end{enumerate}

\subsection{Maximal number of lines}\text{}
\label{ss:maximalsegre}

As a corollary of Theorem \ref{teoremone} and the preceding discussion of cases, we get the following maximality result:

\begin{proposition}\label{prop:aziz}
The maximal numbers of lines on $\cS$ are:
\begin{itemize}
\item $N_d=3d^2$ for $d\geq 3$, $d\not=4,6,8,12,20$;

\item $N_4=64$, $N_6=180$, $N_8=256$, $N_{12}=864$, $N_{20}=1600$.
\end{itemize}
\end{proposition}

\begin{proof}
Looking up at the discussion above, it appears that $\alpha_d=2d$ is maximal
when the group of automorphisms can not be a group $\cT$, $\cO$ or $\cI$ and
that $\alpha_4=12$, $\alpha_6=\alpha_8=24$ and $\alpha_{12}=\alpha_{20}=60$ are
maximal. For other values of $d$, if the automorphism group is $\cT$, resp.
$\cO$, resp. $\cI$ then the number of lines is:
$$
d^2+12d,\quad\text{resp. }d^2+24d,\quad\text{resp. }d^2+60d
$$
and these numbers are bigger than  $3d^2$ only if
$$
d<6,\quad\text{resp. }d<12,\quad\text{resp. }d<30.
$$
So it just remains to check that the degree $d=10$ is not possible for $\cO$ and
$\cI$ and that the degrees $d=14,16,18,22,24,26,28$ are not possible for
$\cI$, that is we cannot decompose such a $d$ as a sum of lengths of orbits for the
groups $\cO$ or $\cI$. This is clear with the restrictions on the numbers of orbits
of each type.
\end{proof}

\begin{remark}
Although this result was expected, one has to pass through the study of \S\ref{ss:possiblesegre} to prove it.
\end{remark}

\subsection{Examples}\text{}
\label{ss:examplessegre}

\begin{enumerate}
\item For $d$ generic, the Fermat surface $F(x,y,z,t)=(x^d-y^d)-(z^d-t^d)$
gives the best example for surfaces of the kind $\phi(x,y)-\psi(z,t)$.

\item For $d=4$, $\Gamma_4\in\{\emptyset,D_2,D_4,\cT\}$ so the possible numbers
of lines for such surfaces are: $16,32,48,64$. This agrees with Segre's result
and $64$ is the maximal possible number of lines on a quartic surface.

\item For $d=5$, $\Gamma_5\in\{\emptyset,\{\id\},C_4,D_3,D_5\}$ so the possible
numbers of lines for such surfaces are: $25,30,45,55,75$. The general bound of
Segre gives $147$.

\item For $d=6$, $\Gamma_6\in\{\emptyset,\{\id\},C_2,D_2,D_3,D_6,\cO\}$ so the
possible numbers of lines for such surfaces are: $36,42,48,60,72,108,180$. The
general bound of Segre gives $240$.

\item The discussion of \S\ref{ss:possiblesegre} gives explicit constructions
of surfaces of each group $\Gamma_d$. For the groups $\cO$ and $\cI$, see also
Section \ref{ch:symmetric}.
\end{enumerate}

\subsection{Real lines}\text{}

It is an interesting problem to find surfaces of any degree $d$ with as many
real lines as possible. For surfaces of the kind $\phi(x,y)-\phi(z,t)=0$, if
the zeros of $\phi$ are all real, one gets already $d^2$ real lines (see proof
of Theorem \ref{teoremone}). Then, for each isomophism in the group $\Gamma_d$
represented by a real matrix, one gets one more real line if $d$ is odd and two
more real lines if $d$ is even.

\section{Surfaces of the kind $t^d=f(x,y,z)$}
\label{ch:geemen}

We consider smooth surfaces of degree $d\geq 3$ given as covering of $\IP_2$
ramified along a plane curve. Let $\cC:f(x,y,z)=0$ be a plane curve defined by
a homogeneous polynomial $f$ of degree $d$ and consider the surface $\cS$ in
$\IP_3$ given by the equation:
$$
F(x,y,z,t):=t^d-f(x,y,z).
$$
Note that the surface $\cS$ is smooth if and only if the curve $\cC$ is.

Set $p=[0:0:0:1]\in \IP_3$. The projection:
$$(\IP_3-\{p\})\rightarrow
\IP_2, \,[x:y:z:t]\mapsto [x:y:z]
$$
induces a $d$-covering $\pi:\cS\rightarrow \IP_2$ ramified along the curve
$\cC$.

Recall that a point $x\in \cC$ is a \emph{$d$-point} (or \emph{total inflection point}) if  the intersection multiplicy of $\cC$ and its tangent line at $x$ is equal to $d$.
\begin{proposition}\label{cov}\text{}
\begin{enumerate}
\item Suppose $L$ is a line contained in $\cS$. Then $\pi(L)$ is a line.

\item \label{splititem2} Let $x\in \cC$ and $L$ the tangent at $\cC$ in $x$, then
the preimage $\pi^{-1}(L)$ consists in $d$ different lines contained in $\cS$ if and only if $x$ is a $d$-point.

\item Let $L$ be a line in $\IP_2$. Then $\pi^{-1}(L)$ contains a line if and
only if $L$ is tangent to $\cC$ at a $d$-point.
\end{enumerate}
\end{proposition}

\begin{proof}\text{}
\begin{enumerate}
\item It is clear from the definition of the projection $\pi$.

\item Assume $x$ is a $d$-point. Let $\Delta$ be a line of equation $\delta$
intersecting $L$ at $x$. Then $d\cdot(\Delta\cdot L)=(\cC\cdot L)$ so after
restriction to $L$ one has up to a scalar factor
$f_{{\mid}_L}=\delta_{{\mid}_L}^d$ showing that the covering restricted to $L$
is trivial and $\pi^{-1}(L)$ consists in the $d$ lines $t-\xi^i\delta_L=0$,
$i=1,\ldots,d$ where $\xi$ is a primitive $d$-th root of the unit. Conversely,
if the covering splits, there exists a section $\gamma\in H^0(L,\cO_{L}(1))$
such that $\gamma^d=f_{{\mid}_L}\in H^0(L,\cO_{L}(d))$ so $L$ intersects $\cC$
at $x$ with multiplicity $d$.

\item If $L$ is the tangent to $\cC$ at a $d$-point the assertion follows from (\ref{splititem2}). Assume now that $\pi^{-1}(L)$ contains a line. Let $L$ be given by a linear function $z=l(x,y)$. Then the equation of $\pi^{-1}(L)$ is $t^d-f(x,y,l(x,y))=0$. Since it contains a line the equation splits as
\begin{eqnarray*}
t^d-f(x,y,l(x,y))=(t-w(x,y))F_{d-1}(t,x,y)
\end{eqnarray*}
where $w(x,y)$ is a linear form. By comparing the coefficients in $t$ one obtains $f(x,y,l(x,y))=w(x,y)^d$ hence the preimage consists in the $d$ lines:
$$
t^d-f(x,y,l(x,y))=\prod_{i=0}^{d-1}(t-\xi^iw(x,y))
$$
where $\xi$ is a primitive $d$-th root of the unit. This means that the covering is trivial over $L$ so by (\ref{splititem2}) $x$ is a $d$-point.
\end{enumerate}

\end{proof}

We deduce the number of lines contained in such surfaces:

\begin{proposition}\label{prop:covering}
Let $\cC:f(x,y,z)=0$ be a smooth plane curve of degree $d$ with $\beta$ total inflection points. Let $\cS$ the surface in
$\IP_3$ given by the equation:
$$
F(x,y,z,t):=t^d-f(x,y,z).
$$
Then $\cS$ contains exactly $\beta\cdot d$ lines. In particular, it contains no more than $3d^2$ lines.
\end{proposition}

\begin{proof}
The first assertion follows directly from the lemma. For the second one, the inflection points are the intersections of $\cC$ with its Hessian curve $\cH$ of degree $3(d-2)$ and at a total inflection point the intersection multiplicity of $\cC$ and $\cH$ is $d-2$, so by Bezout one gets $\beta\leq 3d$.
\end{proof}

\begin{remark}\text{}
\begin{itemize}
\item For $d=3$, it is well-known that each cubic has nine inflection points,
then the induced surface has $3\cdot9=27$ lines.

\item The Fermat curves $x^d+y^d+z^d=0$ have $3d$ total inflection points hence the Fermat surfaces
are examples of surfaces with $3d^2$ lines.
\end{itemize}
\end{remark}

\section{Symmetric surfaces}
\label{ch:symmetric}

We consider surfaces with many symmetries, since one can expect that such
surfaces contain many lines. Indeed, if the surface contains a line then it
contains the whole orbit, and if the symmetry group is big, hopefully this
orbit has big length. To this purpose, we first take $G\subset \PGL(4,\IC)$ be a
finite group of linear transformations acting on $\IP_3$ and construct smooth
$G$-invariant surfaces.

\subsection{Surfaces with cyclic symmetries}\text{}

Denote by $\cU_d$ the group of $d$-th roots of the unit. The group
$\cU_d\times \cU_d$ acts on $\IC[x,y,z,t]$ by $\diag(\xi,\xi,\mu,\mu)$ for
$(\xi,\mu)\in \cU_d\times \cU_d$. The graded space of invariant
polynomials decomposes as:
$$
\IC[x,y,z,t]^{\cU_d\times \cU_d}\cong \IC[x,y]^{\cU_d}\otimes \IC[z,t]^{\cU_d}.
$$
Since $\IC[x,y]^{\cU_d}_k={0}$ for $d\nmid k$ and $\IC[x,y]^{\cU_d}_k=\IC[x,y]_k$
otherwise, all invariant polynomials of degree $d$ for the action of $\cU_d\times
\cU_d$ are of the kind $\phi(x,y)-\psi(z,t)$ for $\phi$ and $\psi$ homogeneous
polynomials of degree $d$. These surfaces were studied in Section
\ref{ch:segre}.

\subsection{Surfaces with polyhedral symmetries}\text{}

We consider again surfaces of the kind $\phi(x,y)=\phi(z,t)$: we studied such
surfaces and their configuration of lines in Section \ref{ch:segre}. We adopt
here a different point of view. Let $\Gamma$ be the group of isomorphisms of
$\IP_1$ permuting the zeros of $\phi$ in $\IP_1$. Then $\phi$ is a projective
invariant for the action of $\Gamma$ on $\IC^2$, \ie
$\phi(g(x,y))=\lambda_g\phi(x,y)$ for $g\in \Gamma$ and $\lambda_g\in \IC^*$.
This implies that the surface $F(x,y,z,t)=\phi(x,y)-\phi(z,t)$ is invariant for
the diagonal action of $\Gamma$ given by $g(x,y,z,t)=(g(x,y),g(z,t))$. Its
number of lines is given by Theorem \ref{teoremone}.

By using this observation, we can find easily equations for surfaces of this
kind with the symmetries of the groups $\cT,\cO,\cI$. The projective invariants
are computed for example in Klein \cite[I.2,\S 11-12-13]{Klein}:

\begin{enumerate}
\item A surface of degree six with octahedral symmetries and $180$ lines:
$$
\phi(x,y)=xy(x^4-y^4).
$$

\item A surface of degree eight with octahedral symmetries and $256$ lines:
$$
\phi(x,y)=x^8+14x^4y^4+y^8.
$$

\item A surface of degree twelve with octahedral symmetries and $432$ lines:
$$
\phi(x,y)=x^{12}-33x^8y^4-33x^4y^8+y^{12}.
$$

\item A surface of degree twelve with icosahedral symmetries and $864$ lines:
$$
\phi(x,y)=xy(x^{10}+11x^5y^5-y^{10}).
$$

\item A surface of degree $20$ with icosahedral symmetries and $1600$ lines:
$$
\phi(x,y)=-(x^{20}+y^{20})+228(x^{15}y^5-x^5y^{15})-494x^{10}y^{10}.
$$

\item A surface of degree $30$ with icosahedral symmetries and $2700$ lines:
$$
\phi(x,y)=(x^{30}+y^{30})+522(x^{25}y^5-x^5y^{25})-10005(x^{20}y^{10}+x^{10}y^{20}).
$$
\end{enumerate}

\subsection{Surfaces with bipolyhedral symmetries}\text{}

First recall the construction of the bipolyhedral groups. Start from the exact
sequence:
$$
0\longrightarrow \{\pm 1\} \longrightarrow \SU(2)
\overset{\phi}{\longrightarrow} \SO(3,\IR)\longrightarrow 0.
$$
For any polyhedral group $G\subset\SO(3,\IR)$, the inverse image
$\widetilde{G}:=\phi^{-1}G$ is called a \emph{binary polyhedral group}. Now
consider the exact sequence:
$$
0\longrightarrow \{\pm 1\} \longrightarrow \SU(2)\times\SU(2)
\overset{\sigma}{\longrightarrow} \SO(4,\IR)\longrightarrow 0.
$$
For $\widetilde{G}$ a binary polyhedral group, the direct image
$\sigma(\widetilde{G}\times\widetilde{G})\subset \SO(4,\IR)$ is called a
\emph{bipolyhedral group}. We shall make use of the following particular
groups:
\begin{itemize}
\item $G_6=\sigma(\widetilde{\cT}\times\widetilde{\cT})$ of order $288$;

\item $G_8=\sigma(\widetilde{\cO}\times\widetilde{\cO})$ of order $1152$;

\item $G_{12}=\sigma(\widetilde{\cI}\times\widetilde{\cI})$ of order $7200$.
\end{itemize}

\bigskip

The polynomial invariants of these groups were studied by Sarti in
\cite{Sarti}. First note that the quadratic form: $Q:=x^2+y^2+z^2+t^2$ is an
invariant of the action of these groups.

\begin{theorem}[Sarti {\cite[\S 4]{Sarti}}]
For $d=6,8,12$, there is a one-dimensional family of $G_d$-invariant surfaces
of degree $d$. The equation of the family is $S_d+\lambda Q^{d/2}=0$. The base
locus of the family consists in $2d$ lines, $d$ in each ruling of $Q$. The
general member of each family is smooth and there are exactly five singular
surfaces in each family.
\end{theorem}

From this theorem immediately follows that each member of the family contains
at least $2d$ lines.

$\bullet$ {\bf The group $G_8$.} Denote by $\cS_8$ the surface $S_8=0$ where:
\begin{align*}
S_8=&x^8+y^8+z^8+t^8+168x^2y^2z^2t^2\\
&+14(x^4y^4+x^4z^4+x^4t^4+y^4z^4+y^4t^4+z^4t^4).
\end{align*}
\begin{proposition}\label{prop:record}
The surface $\cS_8$ contains exactly $352$ lines.
\end{proposition}

\begin{proof}
The proof goes as follows: first we introduce Pl\"ucker coordinates for the lines
in $\IP_3$, then we compute explicitly all the lines contained in the surface.

$\bullet${\it Pl\"ucker coordinates. } Let $\IG(1,3)$ be the Grassmannian of lines in $\IP_3$, or equivalently of
$2$-planes in $\IC^4$. Such a line $L$ is given by a rank-two matrix:
$$
\left(
\begin{matrix}
a & e \\
b & f \\
c & g \\
d & h
\end{matrix}
\right).
$$
The $2$-minors (\emph{Pl{\"u}cker coordinates}):
\begin{align*}
p_{12}&:=af-be & p_{13}&:=ag-ce & p_{14}&:=ah-de\\
p_{23}&:=bg-cf & p_{24}&:=bh-df & p_{34}&:=ch-dg
\end{align*}
are not simultaneously zero, and induce a regular map $\IG(1,3)\longrightarrow
\IP_5$. This map is injective, and its image is the hypersurface
$p_{12}p_{34}-p_{13}p_{24}+p_{14}p_{23}=0$. In order to list once all lines
with these coordinates, we inverse the Pl{\"u}cker embedding in the \emph{Pl\"ucker
stratification}:
$$
\begin{array}{||c||c||c||}
\hline\hline
(1) & (2) & (3) \\
\hline
p_{12}=1 & p_{12}=0, p_{13}=1 & p_{12}=0, p_{13}=0, p_{14}=1\\\hline
\left(
\begin{matrix}
1 & 0 \\
0 & 1 \\
-p_{23} & p_{13} \\
-p_{24} & p_{14}
\end{matrix}
\right)&
\left(
\begin{matrix}
1 & 0 \\
p_{23} & 0 \\
0 & 1 \\
-p_{34} & p_{14}
\end{matrix}
\right) &
\left(
\begin{matrix}
1 & 0 \\
p_{24} & 0 \\
p_{34} & 0 \\
0 & 1
\end{matrix}
\right)\\\hline\hline
(4) & (5) & (6) \\\hline
p_{12}=0, p_{13}=0 & p_{12}=0, p_{13}=0, p_{14}=0 & p_{12}=0, p_{13}=0, p_{14}=0 \\
 p_{14}=0, p_{23}=1 & p_{23}=0, p_{24}=1 & p_{23}=0, p_{24}=0,
p_{34}=1\\\hline \left(
\begin{matrix}
0 & 0 \\
1 & 0 \\
0 & 1 \\
-p_{34} & p_{24}
\end{matrix}
\right)&
\left(
\begin{matrix}
0 & 0 \\
1 & 0 \\
p_{34} & 0 \\
0 & 1
\end{matrix}
\right) &
\left(
\begin{matrix}
0 & 0 \\
0 & 0 \\
1 & 0 \\
0 & 1
\end{matrix}
\right)\\
\hline\hline
\end{array}
$$

$\bullet$ {\it Counting the lines. } The line $L$ is contained in the surface
$\cS_8$ if and only if the function $(u,v)\mapsto S_8(ua+ve,ub+vf,uc+vg,ud+vh)$
is identically zero, or equivalently if all coefficients of this polynomial in
$u,v$ are zero. The conditions for the line to be contained in the surface is
then given by a set of polynomial equations in $a,b,c,d,e,f,g,h$. In order to
count the lines, we restrict the equations to each Pl{\"u}cker stratum and compute
the solutions (this computation is not difficult if left to {\sc singular}
\cite{Sing}).

\begin{enumerate}
\item  {\it The stratum $p_{12}=1$.} Set $p_{23}=c$, $p_{24}=d$, $p_{13}=g$, $p_{14}=h$. The equations for such a line to be contained in the surface are:
\begin{align*}
c^7g+d^7h+7c^3g+7d^3h+7c^4d^3h+7c^3gd^4 &=0\\
c^6g^2+d^6h^2+3c^4d^2h^2+8c^3gd^3h+3c^2g^2d^4&\\
+6c^2d^2+3c^2g^2+3d^2h^2&=0\\
c^5g^3+d^5h^3+c^4dh^3+cg^3d^4+6c^3gd^2h^2&\\
+6c^2g^2d^3h+cg^3+dh^3+6c^2dh+6cgd^2&=0\\
\end{align*}
\begin{align*}
1+g^4+5c^4g^4+5d^4h^4+c^4+d^4+c^4h^4&\\
+g^4d^4+16c^3gdh^3+36c^2g^2d^2h^2+16cg^3d^3h&\\
+h^4+12c^2h^2+12g^2d^2+48cgdh&=0\\
c^3g+d^3h+c^3gh^4+c^3g^5+d^3h^5+6c^2g^2dh^3&\\
+g^4d^3h+6cg^3d^2h^2+6cgh^2+6g^2dh&=0\\
3c^2g^2+3d^2h^2+3c^2g^2h^4+3g^4d^2h^2+c^2g^6&\\
+d^2h^6+8cg^3dh^3+6g^2h^2&=0\\
cg^7+dh^7+7cg^3+7dh^3+7cg^3h^4+7g^4dh^3&=0\\
1+g^8+h^8+14g^4+14h^4+14g^4h^4&=0
\end{align*}

After simplification of the ideal with {\sc singular}
(that we do not reproduce here), the solutions give $320$ lines of the kind $z=cx+gy, t=dx+hy$.
\item {\it The stratum $p_{12}=0$, $p_{13}=1$.} Set $p_{23}=b$, $-p_{34}=d$, $p_{14}=h$. The equations for such a line to be contained in the surface are (after simplification):
\begin{align*}
d&=0\\
b^4h^2-b^2h^4-b^2+h^2&=0\\
b^6-h^6+13b^2-13h^2&=0\\
h^8+14h^4+1&=0\\
b^2h^6+b^4+13b^2h^2+1&=0
\end{align*}
The solutions give $32$ lines of the kind $y=bz, t=hx$, since there are eight possible values for $h$, and for each of them there are four values of $b$.
\end{enumerate}
An easy computation shows that the other strata contain no line, so there are exactly $352$ lines on the surface.
\end{proof}

\begin{remark}
To our knowledge, this is the best example so far of an octic surface with many
lines. This improves widely the bound $256$ of Caporaso-Harris-Mazur
\cite{CHM}.
\end{remark}

$\bullet$ {\bf The group $G_6$.} We take:
$$
S_6=x^6+y^6+z^6+t^6+15(x^2y^2z^2+x^2y^2t^2+x^2z^2t^2+y^2z^2t^2).
$$
\begin{proposition}
The surface $8S_6-5Q^3=0$ contains exactly $132$ lines.
\end{proposition}

There are surfaces with more lines (see \S\ref{ss:examplessegre}), but this
shows the existence of a surface with $132$ lines. This result can be shown in
a similar way as in the $G_8$ case.

\section{A uniform bound}
\label{s:record}

As we mentioned before, the uniform bound $(d-2)(11d-6)$ of Segre is too big
already in degree four. We propose here another lower uniform bound, which
interpolates all maximal numbers of lines known so far, including the octic of
Section \ref{ch:symmetric}. Although there is no reason for this bound to be
maximal, it seems reasonable to expect that an effective construction of a
surface with this number of lines is possible in all degrees.

Let $\cS$ be a smooth surface of degree $d\geq 3$ and $C$ a line contained in
$\cS$. Let $|H|$ be the linear system of planes $H$ passing through $C$. Then
$H\cap \cS=C\cup \Gamma$ where $\Gamma$ is a curve of degree $d-1$. The system
$|\Gamma|$ is described by Segre in \cite{Segre43}: it is base-point free and
any curve $\Gamma$ does not contain $C$ as a component. Then:

\begin{proposition}[Segre {\cite{Segre43}}]\label{prop:segretype}
Either each curve $\Gamma$ intersects $C$ in $d-1$ points which are inflections
for $\Gamma$, or the points of $C$ each of which is an inflection for a curve
$\Gamma$ are $8d-14$ in number. In particular, in this case $C$ is met by no
more than $8d-14$ lines lying on $\cS$.
\end{proposition}

Following Segre, $C$ is called a line of the \emph{second kind} if it
intersects each $\Gamma$ in $d-1$ inflections. A generalization of Segre's
argument in \cite[\S 9]{Segre43} gives the following result:
\begin{proposition}
Assume that $\cS$ contains $d$ coplanar lines, none of them of the second kind. Then $\cS$ contains at most $d(7d-12)$ lines.
\end{proposition}

\begin{proof}
Let $P$ be the plane containing these $d$ distinct lines. Then they are the
complete intersection of $P$ with $\cS$. Hence each other line on $\cS$ must
intersect $P$ in some of the lines. By Proposition \ref{prop:segretype}, each
of the $d$ lines in the plane meets at most $8d-14$ lines, so $8d-14-(d-1)$
lines not on the plane. The total number of lines is at most:
$$
d+d(7d-13)=d(7d-12).
$$
\end{proof}

This bound takes the following values:
$$
\begin{array}{|c||c|c|c|c|c|c|c|c|c|c|}
\hline d & 4& 5& 6& 7& 8& 9& 10& 11& 12& 20\\
\hline\hline 7d^2-12d & 64& 115& 180 & 259& 352& 459& 580& 715& 864& 2560\\
\hline
\end{array}
$$
Note that this bound matches perfectly with the maximal known examples in
degrees $4,6,8,12$.

\section{Number of rational points on a plane curve}
\label{s:aziz}

We give an application of our results to the \emph{universal bound conjecture}, following Caporaso-Harris-Mazur \cite{CHM}:

{\bf Universal bound conjecture. } {\it Let $g\geq 2$ be an integer. There
exists a number $N(g)$ such that for any number field $K$ there are only
finitely many smooth curves of genus $g$ defined over $K$ with more than $N(g)$
$K$-rational points.}

As mentioned in \loccit an interesting way to find a lower bound of $N(g)$, or
of the limit:
$$
\overline{N}:=\limsup_{g\rightarrow \infty}\frac{N(g)}{g}
$$
is to consider plane sections of surfaces with many lines. Indeed, over the
common field $K$ of definition of the surface and its lines, a generic plane
section is a curve containing at least as many $K$-rational points as the
number of lines. In particular, they show that $N(21)\geq 256$. Since we obtain
an octic surface with $352$ lines and a generic plane section of this surface
is a smooth curve of genus $21$, we get:
\begin{corollary}
$N(21)\geq 352$.
\end{corollary}

As we remarked in Section \ref{s:record}, it seems to be possible to construct
surfaces with $d(7d-12)$ lines. This would improve the lower bound of $N(g)$
for many $g$'s. In particular, this would improve the known estimate
$\overline{N}\geq 8$ to $\overline{N}\geq 14$.

\section{Sequences of skew-lines}
\label{s:skew}

A natural question related to the number of lines on a surface is the study of
maximal sequences of pairwise disjoint lines on a smooth surface in $\IP_3$. We
recall the bound of Miyaoka and give some examples.

\subsection{Upper bound for skew-lines}\text{}
\label{ss:upperskew}

The best upper bound known so far for the maximal length of a sequence of
disjoint lines on a smooth surface of degree $d\geq 4$ in $\IP_3$ is given by
Miyaoka:

\begin{theorem}[Miyaoka {\cite[\S 2.2]{Mi}}]
The maximal length of a sequence of skew-lines is $2d(d-2)$ for $d\geq 4$.
\end{theorem}

For $d=3$, each cubic surface contains a maximal sequence of $6$ skew lines.
This comes from the study of the configuration of the $27$ lines (see for
example \cite[Theorem V.4.9]{H} and references therein). For $d=4$, Kummer
surfaces contain a maximal sequence of $16$ skew lines (see for example
\cite{N} and references therein) so the bound is optimal.

But for $d\geq 5$, it is not known if it is sharp.

\subsection{On Miyaoka's bound}\text{}
\label{s:miyaokaproof}

We give a quick sketch of the argument of Miyaoka for the bound on the number
of skew lines, following \cite[\S 2 Examples 2.1,2.2]{Mi}.

Let $X$ be a smooth surface of degree $d\geq 4$ in $\IP_3$. Assume $X$ contains
$r$ disjoint lines $D_1,\ldots,D_r$. By adjunction formula, they have
self-intersection $-n=-(d-2)$. By contracting these lines one gets a surface
$Y$ with $r$ isolated singular points which locally look like the quotient of
$\IC^2$ by a finite group of order $n$.

Write $K_X+\sum\limits_{i=1}^rD_i=P+N'$ with:
\begin{equation*}
P:=K_X+\sum\limits_{i=1}^r \frac{n-2}{n}D_i\quad\text{and}\quad
N':=\sum\limits_{i=1}^r \frac{n-2}{n}D_i.
\end{equation*}
This provides a Zariski decomposition in $\Pic(X)\otimes \IQ$ of
$K_X+\sum\limits_{i=1}^rD_i$.

Set $\nu:=2-1/n$, by \cite[Theorem 1.1]{Mi}, one has the inequality:
$$
r\nu\leq c_2(X)-\frac{1}{3}P^2.
$$
Using that $c_2(X)=d(d^2-4d+6)$ and $K_X^2=d(d-4)^2$ one gets $r\leq 2d(d-2)$.

\subsection{Examples}\text{}
\label{s:geemenrams}

In \cite{Rams05}, Rams considers the surfaces
$x^{d-1}y+y^{d-1}z+z^{d-1}t+t^{d-1}x=0$ and proves that they contain a family
of $d(d-2)+2$ skew-lines for any $d$. In \cite[Example 2.3]{Rams02}, he also
gives an example of a surface of degree five containing a sequence of $19$
skew-lines. We generalize his result, improving the number of skew-lines to
$d(d-2)+4$ in the case $d\geq 7$ and $\gcd (d,d-2)=1$.

Consider the surface $\cR_d:x^{d-1}y+xy^{d-1}+z^{d-1}t+zt^{d-1}=0$. By our
study in Section \ref{s:segregeneral}, this surface contains exactly $3d^2-4d$
lines if $d\neq 6$ and $180$ lines for $d=6$. We prove:

\begin{proposition}\label{essla1}
The surface $\cR_d$ with $\gcd (d,d-2)=1$ contains a sequence of
$d(d-2)+4$ disjoint lines.
\end{proposition}

\begin{proof}
Denote by $\epsilon$, $\gamma$ the primitive roots of the unit of degrees $d-2$
and $d$, and let $\eta:=\epsilon^l\gamma^s$, with $0\leq l\leq d-3$, $0\leq
s\leq d-1$. Since $\gcd (d,d-2)=1$ we have $d(d-2)$ such $\eta$. Now consider
the points
\begin{eqnarray*}
(0:1:0:-\eta^{d-1}), (-\eta:0:1:0)
\end{eqnarray*}
then the line through the two points is
\begin{eqnarray*}
C_{l,s}: (-\eta\lambda:\mu:\lambda:-\eta^{d-1}\mu)
\end{eqnarray*}
An easy computation shows that these lines are contained in $\cR_d$ and are
$d(d-2)$. This form a set of $d(d-2)+4$ skew lines together with the lines
\begin{eqnarray*}
\{x=0, z+\epsilon t=0\},\,\{y=0,z+t=0\},\\
\{z=0,x+\epsilon y=0\},\,\{t=0,x+y=0\}.
\end{eqnarray*}
\end{proof}

\nocite{*}
\bibliographystyle{amsplain}

\bibliography{BoissiereSarti}

\end{document}